\documentclass{scrartcl}

\usepackage{cmap}
\usepackage[utf8]{inputenc}
\usepackage[T1]{fontenc}
\usepackage{microtype}

\usepackage{amsmath}
\usepackage{amssymb}
\usepackage{mathrsfs} 
\usepackage{bm}
\usepackage{csquotes}
\usepackage{enumitem}

\usepackage{multirow}
\usepackage{makecell}
\usepackage{pifont}
\usepackage[table]{xcolor}

\usepackage{orcidlink}

\newcommand{\cmark}{\ding{51}}%
\newcommand{\xmark}{\ding{55}}%
\definecolor{light gray}{gray}{0.90}

\usepackage{tikz}
\usetikzlibrary{automata, matrix, positioning, calc, decorations.pathreplacing, shapes.geometric}
\newlength{\edgelength}
\newcommand{\trans}[4]{%
  \begin{tikzpicture}[auto, shorten >=1pt, >=latex, baseline=(l.base), inner sep=0pt, outer xsep=0.3333em]
    \node (l) {\ensuremath{#1}};%
    \setlength{\edgelength}{\widthof{\scriptsize\ensuremath{#2/#3}}+0.5cm}%
    \node[base right=\edgelength of l] (r) {\ensuremath{#4}};%
    \path[->] (l.mid east) edge node[inner sep=0pt] {\scriptsize\ensuremath{#2/#3}} (r.mid west);%
  \end{tikzpicture}%
}

\usepackage{calc}
\usepackage{tabularx}
\newcommand{\problem}[3][]{%
  \par\vspace{0.125cm plus 0.1cm minus 0.05cm}\begin{tabularx}{\linewidth-2\parindent}{@{}lX}%
    \if\relax\detokenize{#1}\relax%
    \else%
    \textnormal{\textbf{Constant:}}&#1\\%
    \fi%
    \textnormal{\textbf{Input:}}&#2\\%
    \textnormal{\textbf{Question:}}&#3\\%
  \end{tabularx}\vspace{0.125cm plus 0.1cm minus 0.05cm}\par%
}

\usepackage{xspace}

\newcommand*{\SAut}{$\mathscr{S}$\kern-0.4ex-\allowbreak{}au\-to\-ma\-ton\xspace}
\newcommand*{\SAuta}{$\mathscr{S}$\kern-0.4ex-\allowbreak{}au\-to\-ma\-ta\xspace}
\newcommand*{\SIAut}{$\inverse{\mathscr{S}}$\kern-0.4ex-\allowbreak{}au\-to\-ma\-ton\xspace}
\newcommand*{\SIAuta}{$\inverse{\mathscr{S}}$\kern-0.4ex-\allowbreak{}au\-to\-ma\-ta\xspace}
\newcommand*{\GAut}{$\mathscr{G}$\kern-0.2ex-\allowbreak{}au\-to\-ma\-ton\xspace}
\newcommand*{\GAuta}{$\mathscr{G}$\kern-0.2ex-\allowbreak{}au\-to\-ma\-ta\xspace}

\usepackage[affil-it]{authblk}

\author{Emanuele~Rodaro}
\affil{Department of Mathematics\\
  Politecnico di Milano\\
  Piazza Leonardo da Vinci, 32\\
  20133 Milano, Italy}
\author{Jan~Philipp~Wächter \orcidlink{0000-0002-7801-6569}\thanks{The second author is supported by
    CMUP, which is financed by national funds through FCT – Fundação para a Ciência e
    Tecnologia, I.P., under the project with reference UIDB/00144/2020.}}
\affil{Centro de Matemática da Universidade do Porto (CMUP)\\
  Departamento de Matemática, Faculdade de Ciências, Universidade do Porto\\
  Rua do Campo Alegre s/n\\
  4169–007 Porto, Portugal
}

\usepackage{hyperref}
\hypersetup{
  colorlinks,
  linkcolor={red!50!black},
  citecolor={blue!50!black},
  urlcolor={blue!80!black}
}

\title{The Self-Similarity of\\Free Semigroups and Groups}

\begin{document}
  \maketitle
  
  \begin{abstract}
    We give a survey on results regarding self-similar and automaton presentations of free groups and semigroups and related products. Furthermore, we discuss open problems and results with respect to algebraic decision problems in this area.
  \end{abstract}
  
  \begin{section}{Introduction}
    The concept of self-similarity is widespread within mathematics. In this article, we will look at it from an algebraic perspective in terms of groups and semigroups. While this is often done by considering auto- or endomorphisms of infinite regular trees of a certain form (see \cite{nekrashevych2005self} for more background on this), we will choose a different approach based on automaton theory.\footnote{We will only give a very brief introduction. They reader may find more details for example in \cite{waechter2020automaton}, which mostly follows the notation used here.}
    
    \paragraph*{Automata and Runs.}
    In this context, an \emph{automaton}\footnote{From an automaton theoretic point of view, it would be better to actually speak of a letter-to-letter transducer but the simple term \enquote{automaton} is more common in the current context.} is usually a triple $\mathcal{T} = (Q, \Sigma, \delta)$ where $Q$ is a -- typically finite -- set of \emph{states}, $\Sigma$ is an \emph{alphabet} and $\delta \subseteq Q \times \Sigma \times \Sigma \times Q$ is a set of \emph{relations}, for which we will use the more graphical notation $\trans{p}{a}{b}{q}$ to denote a quadruple $(p, a, b, q)$. The intuitive idea is that such a transition of the automaton indicates that, if we start in state $p$ and read an input letter $a$, the automaton outputs the letter $b$ and reaches the state $q$. We may combine multiple transitions, which results in the notion of a run of an automaton. A \emph{run} of the automaton $\mathcal{T} = (Q, \Sigma, \delta)$ is a sequence
    \begin{center}
      \begin{tikzpicture}[auto, >=latex, baseline=(q0.base)]
        \node (q0) {$q_0$};
        \setlength{\edgelength}{\widthof{\scriptsize\ensuremath{a_n/b_n}}+0.5cm}%
        \node[base right=\edgelength of q0] (q1) {$q_1$};
        \node[base right=\edgelength of q1] (dots) {$\dots$};
        \node[base right=\edgelength of dots] (qn) {$q_n$};
        \path[->] (q0.mid east) edge node[inner sep=0pt] {\scriptsize$a_1 / b_1$} (q1.mid west)
                  (q1.mid east) edge node[inner sep=0pt] {\scriptsize$a_2 / b_2$} (dots.mid west)
                  (dots.mid east) edge node[inner sep=0pt] {\scriptsize$a_n / b_n$} (qn.mid west);
      \end{tikzpicture}
    \end{center}
    with $\trans{q_{i - 1}}{a_i}{b_i}{q_i} \in \delta$ for every $1 \leq i \leq n$. It \emph{starts} in $q_0$, \emph{ends} in $q_n$ and its \emph{input} is $a_1 \dots a_n$ and its \emph{output} is $b_1 \dots b_n$. Note that, with our definition of an automaton, the input and output of a run are always of the same length.
    
    We will be dealing both with finite sequences over $\Sigma$ and with finite sequences over $Q$. We call the former \emph{words} and the latter \emph{state sequences} to make a clearer distinction between the two. The set of words (including the empty one) is $\Sigma^*$ and the set of state sequences is $Q^*$. The empty word and the empty state sequence are both denoted by $\varepsilon$.
    
    An automaton $\mathcal{T} = (Q, \Sigma, \delta)$ is \emph{deterministic} if, for every state $p \in Q$ and every letter $a \in \Sigma$, there is \textbf{at most} one transition starting in $p$ with input $a$; i.\,e.\ if
    \[
      \left| \{ \trans{p}{a}{b}{q} \in \delta \mid b \in \Sigma, q \in Q \} \right| \leq 1
    \]
    for all $p \in Q$ and $a \in \Sigma$. Complementary, the automaton $\mathcal{T}$ is \emph{complete} if, for every state $p \in Q$ and every letter $a \in \Sigma$, there is \textbf{at least} one transition starting in $p$ with input $a$; i.\,e.\ if, for all $p \in Q$ and $a \in \Sigma$, we have:
    \[
      \left| \{ \trans{p}{a}{b}{q} \in \delta \mid b \in \Sigma, q \in Q \} \right| \geq 1
    \]
    
    \paragraph*{Automaton Semigroups and Automaton Monoids.}
    In a complete and deterministic automaton $\mathcal{T} = (Q, \Sigma, \delta)$, there is exactly one run starting in $p$ with input $u$ for every state $p \in Q$ and word $u \in \Sigma^*$. We may define $p \circ u$ as the output of this run and $p \cdot u$ as the state it ends in. In this way, every state $q \in Q$ induces a function $\Sigma^* \to \Sigma^*$ mapping $u$ to $p \circ u$ and the closure under composition of these functions gives naturally rise to a semigroup. This is the semigroup $\mathscr{S}(\mathcal{T})$ \emph{generated} by the automaton and any semigroup arising in this way is called a (complete) \emph{automaton semigroup}.
    By adding the identity function on $\Sigma^*$ to the semigroup generated by a (complete and deterministic) automaton $\mathcal{T}$, we obtain a monoid. This is the monoid $\mathscr{M}(\mathcal{T})$ \emph{generated} by $\mathcal{T}$ and (again) every monoid rising in this way is called a (complete) \emph{automaton monoid}.
    
    To extend the notation $p \circ u$ with $u \in \Sigma^*$ to multiple states, we let $p_n \dots p_2 p_1 \circ u = p_n \circ ( \dots p_2 \circ (p_1 \circ u))$ for $p_1, p_2, \dots, p_n \in Q$ and $\varepsilon \circ u = u$ for the empty state sequence $\varepsilon$ (which corresponds to adding the identity function in the case of the generated monoid).
    
    \paragraph*{Automaton Groups.}
    An automaton $\mathcal{T} = (Q, \Sigma, \delta)$ is \emph{invertible} if, for every state $p \in Q$ and every letter $b \in \Sigma$, there is at most one transition starting in $p$ with output $b$; i.\,e.\ if
    \[
      \left| \{ \trans{p}{a}{b}{q} \in \delta \mid b \in \Sigma, q \in Q \} \right| \leq 1
    \]
    for all $p \in Q$ and $b \in \Sigma$. If the automaton is not only invertible but also complete, the above set will contain exactly one element for every state $p \in Q$ and (output) letter $b \in \Sigma$ by reasons of cardinality. From this, it is not difficult to see that the functions induced by the states of an invertible, deterministic and complete automaton are bijections $\Sigma^* \to \Sigma^*$. Taking the closure of these functions and their inverses under composition yields a group, which is the group $\mathscr{G}(\mathcal{T})$ \emph{generated} by $\mathcal{T}$ and any such group is called an \emph{automaton group}.
    
    \paragraph*{Partial Automata and Inverse Automaton Semigroups.}
    If we drop the requirement of being complete and only consider deterministic automata $\mathcal{T} = (Q, \Sigma, \delta)$, we still obtain that, for every state $p \in Q$ and every input $u \in \Sigma^*$, there is \textbf{at most} one transition starting in $p$ with input $u$. This means that the functions induced by the states are now \textbf{partial} functions $\Sigma^* \to \Sigma^*$ (and the same is true for their compositions, of course). We may still define the semigroup or monoid generated by such an automaton though, and obtain the notion of a \emph{partial automaton semigroup} or \emph{monoid}. We will not go into more detail about this concept (as it is not widely studied in the literature) but the reader may find an introduction in \cite{structurePart}.
    
    The advantage of this approach is (not only that it is arguably more natural in the setting of semigroups but also) that it allows to consider deterministic, invertible but possibly non-complete automata. The functions induced by their states are partial injections, which are closely related to \emph{inverse} semigroups.\footnote{We will not elaborate on inverse semigroup but refer the reader to the standard literature, e.\,g.\ \cite{petrich1984, HOWIE}.} This leads to the notion of an \emph{inverse automaton semigroup}. There are some subtle points to consider here,\footnote{In particular, there is a priori a difference between automaton semigroups which happen to be inverse and the semigroups generated by invertible, deterministic automata.} for which refer the reader again to \cite{structurePart}.
    
    \begin{table}\centering
      \renewcommand\cellalign{lc}
      \setlength{\extrarowheight}{0.25ex}
      \begin{tabular}{cl|cccc}
        \multirow[c]{3}{*}{\rotatebox{90}{$\mathcal{T}$ is\dots}} & deterministic & \multicolumn{4}{c}{\cellcolor{light gray}Yes} \\
        & complete & \multicolumn{2}{c}{\cellcolor{light gray}Yes} & \multicolumn{2}{c}{No} \\
        & invertible & \cellcolor{light gray}Yes & No & \cellcolor{light gray}Yes & No \\\hline\hline
        
        \multirow{4}{*}[-1ex]{\rotatebox{90}{$\mathcal{T}$ generates\dots}} & \makecell{partial automaton\\semigroup/monoid} & \cmark & \cmark & \cmark & \cmark \\\cline{2-6}
        & \makecell{(complete) automaton\\semigroup/monoid} & \cmark & \cmark & \xmark & \xmark \\\cline{2-6}
        & \makecell{inverse automaton\\semigroup/monoid} & \cmark & \xmark & \cmark & \xmark \\\cline{2-6}
        & \makecell{automaton group} & \cmark & \xmark & \xmark & \xmark \\
      \end{tabular}\caption{The various algebraic structures generated by different kinds of automata.}\label{tbl:generatedStructures}
    \end{table}
    \paragraph*{Automaton Structures and Their Self-Similar Nature.}
    In order to keep our terminology lightweight, we will use the term \emph{automaton structures} to refer to the various concepts defined above (automaton semigroups, monoids, groups, inverse semigroups; complete and partial; \dots). An overview of these and the properties required in the generating automaton may be found in \autoref{tbl:generatedStructures}.
    
    The way we have defined them, automaton structures consist of length-preserving, prefix-compatible functions. The latter means that $\bm{p} \circ uv = u' v'$ implies $\bm{p} \circ u = u'$ for all state sequences $\bm{p}$ and words $u, u', v, v'$ where $u$ and $u'$ are of the same length. Thus, $\bm{p}$ and $u$ uniquely yield a function $\Sigma^* \to \Sigma^*$ mapping $v$ to $(\bm{p} \cdot u) \circ v = v'$ which is determined by $\bm{p} \circ uv = u'v'$. We say that this is the function induced by $\bm{p}$ \emph{shifted} by $u$. In an automaton structure, this shifted function is again given by a state sequence, which meas that it is itself contained in the automaton structure. This is the motivation for saying that automaton structures are \emph{self-similar}: the shifted functions are contained in the structure and, thus, \emph{similar} to the non-shifted functions.
    
    In fact, the state sequences inducing the function belonging to a state sequence $p_n \dots p_1$ (with $p_1, \dots, p_n \in Q$) shifted by $u$ can be defined inductively by $p_n \dots p_2 p_1 \cdot u = \left( p_n \dots p_2 \cdot \left( p_1 \circ u \right) \right) (p_1 \cdot u)$ (where the definition of $p \cdot u$ for a single state is the one stated above). The reader may verify that this definition indeed satisfies $\bm{p} \circ uv = (\bm{p} \circ u) \left( (\bm{p} \cdot u) \circ v \right)$ for state sequences $\bm{p}$.
    
    \paragraph*{(General) Self-Similar Structures.}
    In an automaton structure, if we shift the function induced by a state, the resulting function will again be induced by a state (of which we only have finitely many). In general, however, self-similarity only means that the shifted version of a semigroup (monoid, group) element is again in the semigroup (monoid, group). To cover this case, we can simply drop the requirement for our automata to only contain finitely many states. Then, we may have a state for every semigroup (monoid, group) element. The structures generated by such infinite-state automata are called \emph{self-similar} and we obtain a self-similar version for each of the automaton structures listed in \autoref{tbl:generatedStructures}: self-similar semigroups/monoids, self-similar inverse semigroups/monoids, their partial counter-parts and self-similar groups.
  \end{section}

  \begin{section}{The History of Free Structures and Self-Similarity}
    There is a long history on the problem of presenting free groups and semigroups in a self-similar way. The first construction for presenting a free, non-abelian group\footnote{The free group in one generator is generated by the \emph{adding machine}; see, e.\,g.\ \cite[Example~0.2.1.4]{waechter2020automaton}.} in a self-similar way or even as an automaton groups seems to be the \emph{Aleshin automaton} \cite{aleshin1983free}. However, it was only proved to generate a free group of rank three much later by Vorobets and Vorobets \cite{vorobets2007free}. Prior to that result, Brunner and Sidki presented the general linear group $\operatorname{GL}_n(\mathbb{Z})$ (and, thus, also free groups of finite or countable rank) as a subgroup of an automaton group \cite{brunner1998generation}. Shortly after, free groups were presented as subgroups of an automaton group over a binary alphabet by Oliinyk and Sushchansky \cite{Oliynyk2000free}. The first proof for an automaton to generate a free group was given by Glasner and Mozes \cite{glasner2005automata}. However, in contrast to Aleshin's automaton, the rank of the generated group is not the same as the number of states (there are twice as many states). Automata generating free groups of odd rank (starting at rank three) where the rank coincides with the number of states were constructed from the Aleshin automaton by Vorobets and Vorobets \cite[Theorem~1.3]{vorobets2010series}. These automata can also be combined to obtain (disconnected) automata for free groups of arbitrary rank (except some small numbers) \cite[Theorem~1.4]{vorobets2010series}. Together with Steinberg, Vorobets and Vorobets later also derived a family of connected automata generating free groups of even rank (where the rank is at least four and coincides with the number of states) \cite{steinberg2011automata}.
    
    Compared to free groups, the situation for free semigroups is much simpler: the free semigroup of rank one is not an automaton semigroup \cite[Proposition~4.3]{cain_1auto} but free semigroups of higher rank are \cite[Proposition~4.1]{cain_1auto}. In fact, the construction used to generate these free semigroups is surprisingly simple. On the other side, the argument used to show that the free semigroup of rank one is not an automaton semigroup has various generalisations (see \cite[Theorem~15]{bc_automaton2}, \cite[Theorem~19]{structurePart} and \cite[Theorem~1.2.1.4]{waechter2020automaton}). However, a generalisation for showing that the semigroup is not self-similar is not obvious and results in this direction do not appear to exist. In the case of monoids, the free monoid of rank one is indeed an automaton monoid (it is generated by the adding machine, just like the free group of rank one) and free monoids of higher rank can be generated in the same way as the corresponding free semigroups. Finally, for inverse semigroups and monoids, there does not seem to exist much research. Based on a presentation of the monogenic free inverse semigroup as a subsemigroup of an inverse automaton semigroup by Oliynyk, Sushchansky and Slupik \cite[Theorem~25]{Oliynyk2010inverse}, D'Angeli and authors presented this semigroup as an inverse automaton semigroups \cite[Example~2]{decidabilityPart} (see also \cite[Example~23]{structurePart}).
    
    Related to the question on how to present free groups and semigroups in a self-similar way, there is another line of research on presenting free products of groups or semigroups self-similarly. The free product of two groups or semigroups $X = \langle P \mid \mathcal{R} \rangle$ and
    $Y = \langle Q \mid \mathcal{S} \rangle$ is the group or semigroup $X \star Y = \langle P \cup Q \mid \mathcal{R}\cup \mathcal{S} \rangle$. Depending on whether we want to take the free product in the category of groups or in the category of semigroups, we understand this as group or semigroup presentation. Note that this distinction is of importance since the free product of two groups in the category of semigroups is in general not a group (since there is no neutral element)! Similarly to the Aleshin automaton above, there is the Bellaterra automaton, which generates the free product of three groups of order two \cite[Theorem~1.10.2]{nekrashevych2005self} -- a result due to Muntyan and Savchuk. This group had previously been presented as a subgroup of an automaton group \cite[p.~323]{Oliynyk2000free}. Note, however, that every free product of finitely many finite groups is a subgroup of an automaton group.\footnote{The authors would like to thank Armin Weiß for pointing this out!} In fact, this is also true for amalgamated products and HNN extensions.\footnote{More on amalgamated products and HNN extension can be found, for example, in \cite{lyndon2001combinatorial}.} This follows from the fact that such products are virtually free.\footnote{They are fundamental groups of finite graphs of finite groups, see \cite{diekert2017contextfree} for an introduction.} Since the free group in two generators (and, thus, every free group of countable rank) is a subgroup of $\operatorname{SL}_2(\mathbb{Z})$ (by a well-known embedding), virtually free groups are, in particular, subgroups of $\operatorname{GL}_n(\mathbb{Z})$ (which can be seen from their induced representation). Therefore, by Brunner and Sidki's result \cite{brunner1998generation}, they are subgroups of an automaton groups. However, there are also explicit (often much simpler) constructions to obtain such groups as subgroups of automaton groups. For example, Gupta, Gupta and Oliynyk gave a construction for free products of finite groups \cite{gupta2007free} and there are further constructions for free products of finite cyclic groups with amalgamation over a common cyclic subgroup \cite{prokhorchuk2021generation, oliynyk2021amalgamated}. In addition, there are also results for amalgamated products and HNN extensions of non-finite groups. For example, Lavrenyuk, Mazorchuk, Oliynyk and Sushchansky presented free products of two infinite cyclic groups with amalgamation over an infinite cyclic subgroup as a subgroup of an automaton group \cite{lavrenyuk2007faithful} and Prokhorchuk considered HNN extensions of certain free abelian groups \cite{prokhorchuk2021finite}. With regard to self-similar presentations, Vorobets and Vorobets generalized the Bellaterra automaton to present any free product of an odd number of groups of order two as an automaton group \cite[Theorem~1.7~(i)]{vorobets2010series}; as with the Aleshin automaton, these automata can be combined to also cover (sufficiently large) even numbers of copies of the group of order two \cite[Theorem~1.7~(ii)]{vorobets2010series}. This was extended to an arbitrary number of copies by Savchuk and Vorobets \cite[Theorem~0.2]{savchuk2011automata} (see also \cite{steinberg2011automata}).
    
    More fundamentally, Fedorova and Oliynyk showed that the free product of finitely many subgroups of an automaton group is again a subgroup of an automaton group \cite{fedorova2017finite}. However, there do not seem to be results on presenting arbitrary free products of self-similar groups in a self-similar way. For semigroups, on the other hand, such results exist.
    Probably the first results in this direction were given by Brough and Cain: first they showed that the free product of two automaton semigroups both containing a left identity is again an automaton semigroup \cite[Theorem~2]{bc_automaton1} and, second, they showed that, if $S$ and $T$ are automaton semigroups, then the semigroup $(S \star T)^1$, which arises by taking the free product of $S$ and $T$ and adjoining a neutral element, is again an automaton semigroup \cite[Theorem~3]{bc_automaton1}. Note, however, that there is again a subtle difference between semigroups and monoids here: if $S$ and $T$ are monoids their free product of monoids is different in general to their free product as semigroups with an additional identity (compare to \cite[Corollary~4]{bc_automaton1})! They could later extend this to show that the free product of two automaton semigroups that either both contain an idempotent or are both homogeneous (with respect to their automaton presentation)\footnote{An automaton semigroup is homogeneous with respect to its automaton presentation if two state sequences can only induce the same function if they are of the same length.} is an automaton semigroup \cite[Theorem~4]{bc_automaton2}.
    The construction from this proof was later adapted in Welker's Bachelor thesis \cite{welker} to relax the hypothesis further. Using a modified construction, one can even show that the free product of any two automaton semigroups such that there is a homomorphism from one to the other is again an automaton semigroup \cite{brough2022automaton}.
  \end{section}

  \begin{section}{Computational Aspects}
    The field of algebraic decision problems goes back to a seminal paper by Max Dehn from 1911 \cite{dehn11}. In this paper, he formulated three problems, which have become known as the three fundamental problems of algorithmic group theory. The first one is the \emph{word problem}, which asks whether a given word over the generators of a group is the neutral element.\footnote{For a discussion of the word problem with respect to automaton structures see \cite{dangeli2017complexity, waechter2022automaton}.} The second one is the \emph{conjugacy problem}: given two group elements, decide whether they are conjugated.\footnote{There is an automaton group with an undecidable conjugacy problem \cite{sunic2012conjugacy}.} The third one, the \emph{isomorphism problem}, however, is of a slightly different nature: it asks whether two given groups are isomorphic. The fundamental difference is that, for this problem, the groups are part of the input, which requires them to be provided in a suitable, finite (or at least recursively enumerable) presentation. Typically, this is done by providing a finite set of generators and a finite set of relations. The Adian-Rabin theorem states that most \enquote{reasonable} properties of thus provided groups are undecidable (see, e.\,g.\ \cite{lyndon2001combinatorial}).
    
    Alternatively, we may also use complete, deterministic and invertible automata to provide groups to algorithms and it turns out that not much is known on algorithmic problems for this kind of presentation. While the isomorphism problem for automaton groups is known to be undecidable (which follows from \cite{sunic2012conjugacy}), the decidability of very similar problems such as the \emph{finiteness problem} (\enquote{is a given group finite?}) remain open \cite[7.2.(b)]{grigorchuk_automata}, moving a more general result such as an analogue of the Adian-Rabin theorem into even further distance. Typically, these kinds of problems have natural generalization to (possibly inverse) semigroup and monoids where some conditions on the generating automaton are relaxed. This usually facilitates the encoding of Turing machines, which makes obtaining results easier. For example, the finiteness problem for automaton semigroups is undecidable \cite{gillibert_finiteness}. While interesting in its own right, the generalization to semigroups is often done in the hope to later extend results to groups.\footnote{For example, the \emph{order problem} -- which we will not discuss further -- was first known to be undecidable for automaton semigroups \cite{gillibert_finiteness} and could later be proved to be also undecidable for automaton groups \cite{gillibert2017automaton}.}
    
    A problem of particular interest in our current context is the \emph{freeness problem} for automaton groups:
    \problem{%
      a complete, invertible and determinstic automaton $\mathcal{T}$
    }{%
      is $\mathscr{G}(\mathcal{T})$ free?
    }\noindent
    Just like with the finiteness problem, the decidability of this problem is still open \cite[section 2, 2.(g)]{AIMproblemList2007} \cite[7.2.(b)]{grigorchuk_automata}. However, there are some partial results also with regard to the very similar problem whether the state set of a given automaton forms a free basis in the generated group:
    \problem{%
      a complete, invertible, deterministic automaton with states $Q$
    }{%
      is $\mathscr{G}(\mathcal{T})$ free with a basis induced by $Q$?
    }\noindent
    It also makes sense to consider a slight variation of this problem where one dedicated state of the automaton induces the identity function and we exclude this state as an element of the basis.
    
    Clearly, $Q$ yields a basis of the generated group if no (non-empty) sequence of states and their inverses induces the identity (and, thus, the neutral element of the group), i.\,e.\ if we have no group relation over $Q$. Together with D'Angeli, the current authors have shown a result which is related to this view: it is not possible to decide whether a given complete, invertible and determinstic automaton admits a (non-empty) sequence over the states (but not their inverses) which induces the identity function    \cite[Theorem~3.7]{decidabilityPart}. This has some interesting consequences as it implies that it is undecidable whether a given automaton semigroup contains a neutral element (see \cite[Subsection 2.2.2]{waechter2020automaton} for a discussion).

    This result already seems very close to showing that the freeness problem for automaton \textbf{semigroups}
    \problem{%
      a complete, determinstic automaton $\mathcal{T}$
    }{%
      is $\mathscr{S}(\mathcal{T})$ free?
    }\noindent
    is undecidable. However, the situation is more complex and this particular approach does not seem to yield the desired result \cite{decidabilityPartErratum}. On the other hand, it is known that the problem whether certain invertible, complete and deterministic automata with only two states generate a free semigroup is decidable \cite{klimann2016automaton}.
    
    In general, algebaric decision problems on freeness seem to be closely related to the classical Post Correspondence Problem.
    For example, the above mention result on finding a sequence of (positive) states inducing the identity is obtained by using a reduction from some variation of this problem.
  \end{section}
  
  \bibliographystyle{plain}
  \bibliography{references}
\end{document}